\documentclass[12pt,a4paper]{amsart}
\usepackage{amssymb,amsmath,amsthm}
\usepackage{graphicx}
\usepackage{hyperref}
 \usepackage[color=green!40]{todonotes}

\makeatletter
\newtheorem*{rep@theorem}{\rep@title}
\newcommand{\newreptheorem}[2]{%
\newenvironment{rep#1}[1]{%
 \def\rep@title{#2~\ref{##1}}%
 \begin{rep@theorem}}%
 {\end{rep@theorem}}}
\makeatother

\theoremstyle{plain}
\newtheorem{theorem}{Theorem}%[section]
\newreptheorem{theorem}{Theorem}

\newtheorem{observation}[theorem]{Observation}

\theoremstyle{definition}

\newtheorem{fact}[theorem]{Fact}

\newcommand\cref[1]{Corollary~\ref{cor:#1}}

\DeclareMathOperator{\mis}{mis}

\title{On the number of maximal independent sets: From Moon-Moser to Hujter-Tuza}

\author{Cory Palmer}
\address{University of Montana}
\email{cory.palmer@umontana.edu}
\thanks{Palmer's research is supported by a grant from the Simons Foundation \#712036}

\author{Bal\'azs Patk\'os}
\address{Alfr\'ed R\'enyi Institute of Mathematics}
\email{patkos@renyi.hu}
\thanks{Patk\'os's research is partially supported by NKFIH grants SNN 129364 and FK 132060}

\date{}
\begin{document}

\maketitle

\begin{abstract}
We connect two classical results in extremal graph theory concerning the number of maximal independent sets. The maximum number $\mis(n)$ of maximal independent sets in an $n$-vertex graph was determined by Moon and Moser. The maximum number $\mis_\bigtriangleup(n)$ of maximal independent sets in an $n$-vertex triangle-free graph was determined by Hujter and Tuza. We determine the maximum number $\mis_t(n)$  of maximal independent sets in an $n$-vertex graph containing no induced triangle matching of size $t+1$. We also reprove a stability result of Kahn and Park on the maximum number $\mis_{\bigtriangleup,t}(n)$ of maximal independent sets in an $n$-vertex triangle-free graphs containing no induced matching of size $t+1$.
\end{abstract}

\section{Introduction}

Let $\mis(G)$ denote the of maximal independent sets in the graph $G$. The classic result determining $\mis(n)$, the maximum of $\mis(G)$ over all graphs on $n$ vertices is:

\begin{theorem}[Miller, Muller \cite{MiMu}, Moon, Moser \cite{MM}]\label{mm}
For any $n\ge 3$, we have 
\begin{eqnarray*}
\mis(n)=\left\{
\begin{array}{cc} %rr ll cr lc
3^{n/3} & \textnormal{if} ~n ~\textnormal{is divisible by 3}, \\
4\cdot 3^{(n-4)/3}& \textnormal{if $n\equiv 1$ (mod 3)}, \\
2\cdot 3^{(n-2)/3}& \textnormal{if $n\equiv 2$ (mod 3)}.
\end{array}
\right.
\end{eqnarray*}
\end{theorem}

For graphs $F,G$ and positive integers $a,b$, we denote by $aF+bG$ the vertex-disjoint union of $a$ copies of $F$ and $b$ copies of $G$. 
Then the constructions giving the lower bound of Theorem~\ref{mm} are $\frac{n}{3}K_3$, $\frac{n-4}{3}K_3+K_4$ or $\frac{n-4}{3}K_3+2K_2$, and $\frac{n-2}{3}K_3+K_2$ in the three respective cases. As these constructions contain many triangles, one can ask the natural question to maximize $\mis(G)$ over all triangle-free graphs. The maximum over all such $n$-vertex graphs, denoted by $\mis_\bigtriangleup(n)$, was determined by Hujter and Tuza~\cite{HT}.

\begin{theorem}[Hujter, Tuza \cite{HT}]\label{ht}
For any $n\ge 4$, we have 
\begin{eqnarray*}
\mis_{\bigtriangleup}(n)=\left\{
\begin{array}{cc} %rr ll cr lc
2^{n/2} & \textnormal{if} ~n ~\textnormal{is even}, \\
5\cdot 2^{(n-5)/2}& \textnormal{if $n$ is odd}.
\end{array}
\right.
\end{eqnarray*}
\end{theorem}

The parameter $\mis(G)$ has been determined for connected graphs (see~\cite{F,Gr}) and for trees (see~\cite{Wi,Sa}). The value of $\mis(n)$ has implications for the runtime of various graph-coloring algorithms (see \cite{W} for several references).

Answering a question of Rabinovich, Kahn and Park~\cite{KP} proved stability versions of both Theorem~\ref{mm} and Theorem~\ref{ht}. An {\it induced triangle matching} is an induced subgraph that is a vertex disjoint union of triangles; its {\it size} is the number of triangles.

\begin{theorem}[Kahn, Park \cite{KP}]\label{kp}
For any $\varepsilon >0$,  there is a $\delta=\delta(\varepsilon) = \Omega(\varepsilon)$ such that for any $n$-vertex graph $G$ that does not contain an induced triangle matching of size $(1-\varepsilon)\frac{n}{3}$, we have $\log \mis(G)<(\frac{1}{3}\log 3-\delta)n$.
\end{theorem}

An {\it induced matching} is an induced subgraph that is a matching; its {\it size} is the number of edges

\begin{theorem}[Kahn, Park \cite{KP}]\label{kp2}
For any $\varepsilon >0$,  there is a $\delta=\delta(\varepsilon) = \Omega(\varepsilon)$ such that for any $n$-vertex triangle-free graph $G$ that does not contain an induced matching of size $(1-\varepsilon)\frac{n}{2}$, we have $\log \mis(G)<(\frac{1}{2}-\delta)n$.
\end{theorem}

Let $\mis_t(n)$ denote the maximum number of maximal independent sets in an $n$-vertex graph that does not contain an induced triangle matching of size $t+1$. With this notation we have $\mis_0(n) = \mis_\bigtriangleup(n)$ and Theorem~\ref{kp} gives $\mis_t(n) < 3^{(1/3-\delta')n}$ when $t+1 = (1-\epsilon)\frac{n}{3}$.

The primary result of this note is the following common generalization of Theorems~\ref{mm} and \ref{ht} which gives a strengthening of Theorem~\ref{kp} as it determines $\mis_t(n)$ for all $n$ and $t\leq n/3$. 
% \corys{should we state more explicitly how this gives a stronger version of Thm 3?}
% \pb{We could continue the sentence as 'by determining the exact value of $\mis_t(n)$ for every pair $n,t$ of integers with $0\le t\le n/3$. Observe that with our notation Theorem \ref{kp} states that if $t\le (1-\varepsilon)\frac{n}{3}$, then $\mis_t(n)\le (3-\delta)^n$.'}

\begin{theorem}\label{main}
For any $0\le t \le n/3$ put $m=n-3t$. Then we have 
\begin{eqnarray*}
\mis_{t}(n)=\left\{
\begin{array}{cc} %rr ll cr lc
3^t\cdot 2^{m/2} & \textnormal{if} ~m ~\textnormal{is even}, \\
3^{t-1}\cdot 2^{(m+3)/2}& \textnormal{if $m$ is odd and $t>0$},
\\
5\cdot 2^{(n-5)/2}& \textnormal{if $m$ is odd and $t=0$}.
\end{array}
\right.
\end{eqnarray*}
\end{theorem}

Constructions showing the lower bounds are $tK_3+\frac{m}{2}K_2$, $(t-1)K_3+\frac{m+3}{2}K_2$, and $C_5+\frac{n-5}{2}K_2$, respectively.

Let $\mis_{\bigtriangleup,t}(n)$ denote the maximum of $\mis(G)$ over all $n$-vertex triangle-free graphs that do not contain an induced matching of size $t+1$. The secondary result of this note is the following short reproof of Theorem~\ref{kp2}.

\begin{theorem}\label{kp2+}
Let $c$ denote the largest real root of the equation $x^6-2x^2-2x-1=0$, $c=1.40759\ldots<\sqrt{2}$. Then $\mis_{\bigtriangleup,t}(n)\le 2^tc^{n-2t}$.
\end{theorem}

\section{Proofs}

In our proofs we shall use an observation due to  Wood~\cite{W}. It follows from the fact that any maximal independent set in $G$ must meet the closed neighborhood $N[v]=N(v) \cup \{v\}$ of any vertex $v$ of $G$.

\begin{observation}[Wood \cite{W}]\label{obs}
For any graph $G$ and vertex $v\in V(G)$ we have 
\[
\mis(G)\le \sum_{w\in {N[v]}}\mis(G\setminus N[w]).
\]
% Thus, if $G$ has minimum degree $d$, then
% \[
% \mis(G) \le d \cdot \mis_t(n-d-1).
% \]
\end{observation}

We begin with some inequalities involving the bound, denoted $g_t(n)$, in Theorem~\ref{main}.
% We introduce the function $g_t(n)$ for all values $0\le t,n$ such that if $t\le n/3$, then $g_t(n)$ equals the bound in Theorem~\ref{main}:

% \begin{eqnarray*}
% g_{t}(n)=\left\{
% \begin{array}{cc} %rr ll cr lc
% 3^{t^*}\cdot 2^{m/2} & \textnormal{if} ~m ~\textnormal{is even}, \\
% 3^{t^*-1}\cdot 2^{(m+3)/2}& \textnormal{if $m$ is odd and $t>0$},
% \\
% 5\cdot 2^{(n-5)/2}& \textnormal{if $m$ is odd and $t=0$},
% \end{array}
% \right.
% \end{eqnarray*}
% where $t^*:=\min\{t,\lfloor n/3\rfloor\}$.

\begin{fact}\label{fact}
For $t>0$, we have $\frac{g_t(n-3)}{g_t(n)}\le 3/8$, $\frac{g_t(n-2)}{g_t(n)}=1/2$, and $\frac{g_t(n-4)}{g_t(n)}=1/4$.
\end{fact}

Observe that if $n$ is odd and $t=0$, then the bounds in Fact~\ref{fact} may not hold and, in particular, Case III of the following argument will not work. Fortunately, we may assume that $t >0$ as the $t=0$ case is exactly Theorem~\ref{ht}.

% then $\frac{g_0(n-3)}{g_0(n)}=2/5$, so $2g_0(n-3)+g_0(n-4)>g_0(n)$. Therefore, Case III of our proof will not work if $t=0$. The proof of Hujter and Tuza will complete our proof. 

In the proof, we will always compare $g_t(n-k)$ to $g_t(n)$, and it might happen that $n-k$ drops below $3t$. In this case, we consider $g_t(n-k)$ to be $g_{\lfloor \frac{n-k}{3}\rfloor}(n-k)$. Fortunately, all inequalities in Fact~\ref{fact} remain true as $g_t(n)$ is non-decreasing in $t$.

\begin{proof}[Proof of Theorem~\ref{main}]
By the discussion above we may assume $t>0$.
We proceed by induction on $m$. Observe that cases $m=0,1,2$ are covered by Theorem~\ref{mm}. %(and also case $m=4$ by the construction $\frac{n-4}{3}K_3+2K_2$). 
Let $G$ be a graph on $n$ vertices not containing an induced triangle matching of size $t+1$. We distinguish cases according to the minimum degree of $G$.

\medskip

\textsc{Case I:} $G$ has a vertex $x$ of degree $1$.

\medskip

Then by applying Observation~\ref{obs} with $v=x$ and Fact~\ref{fact}, we obtain $\mis(G)\le 2\mis_t(n-2) \leq 2g_t(n-2) \leq g_t(n)$.

\medskip

\textsc{Case II:} $G$ has a component $C$ of minimum degree $d \geq 3$.

\medskip

Then by applying Observation~\ref{obs} to any $v\in C$ and Fact~\ref{fact}, we obtain $\mis(G)\le (d+1) \mis_t(n-d-1)\le g_t(n)$.

\medskip

\textsc{Case III:} $G$ has component $C$ with a vertex $x$ of degree $2$ and a vertex of degree at least $3$.

\medskip

We may assume that $x$ is adjacent to a vertex $y$ of degree $d(y) \geq 3$.
 Applying Observation~\ref{obs} with $v=x$ and Fact~\ref{fact}, we obtain $\mis(G)\le 2 \mis_t(n-3)+\mis_t(n-4)\le 2g_t(n-3) + g_t(n-4)\leq g_t(n)$. 

\medskip

\textsc{Case IV:} $G$ is $2$-regular, i.e., a cycle factor.

\medskip

It is not hard to verify (see for example \cite{F}) that $\mis(C_3)=3, \mis(C_4)=2, \mis(C_5)=5$ and $\mis(C_n)=\mis(C_{n-2})+\mis(C_{n-3})$. In particular, if $n\neq 3$, then $\mis(C_n)^{1/n}$ is maximized for $n=5$ with value $5^{1/5}$. Thus, for cycle factors containing at most $t$ triangles, we have $\mis(G)\le 3^t\cdot 5^{(n-3t)/5}\le g_t(n)$.
\end{proof}

Before the proof of Theorem \ref{kp2+}, we gather facts about the bound $h_t(n):=2^t\cdot c^{n-2t}$.

\begin{fact}\label{fact2}
For the largest real root $c=1.40759\ldots<\sqrt{2}$ of $x^6-2x^2-2x-1=0$, we have
\begin{enumerate}
    \item 
    $2+c\le 2c^2$ and so $h_t(n-2)+h_{t-1}(n-3)\le h_t(n)$,
    \item
    for any $d\ge 4$, we have $(d+1)\le c^{d+1}$ and so $(d+1)h_t(n-d-1)\le h_t(n)$,
    \item
    $3c+1\le c^5$ and so $3h_t(n-4)+h_t(n-5)\le h_t(n)$,
    \item
    $2c+1\le c^4$ and so $2h_t(n-3)+h_t(n-4)\le h_t(n)$.
\end{enumerate}
\end{fact}

Just as in Fact \ref{fact}, $n-k$ might drop below $2t$, in which case we consider $h_t(n-k)$ to be $h_{\lfloor \frac{n-k}{2}\rfloor}(n-k)$ and again the inequalities in Fact~\ref{fact2} remain true.

\begin{proof}[Proof of Theorem~\ref{kp2+}]
We use double induction. First on $m:=n-2t$ and then on $t$. Theorem~\ref{ht} yields the statement when $m=0,1$. For any $m$, if $t=0$, then the only graph $G$ on $m-2t=m$ vertices that does not contain an induced matching of size one is the empty graph and in this case $\mis(G)=1\le h_0(m)=c^{m}$. Let $G$ be an $n$-vertex triangle-free graph that contains no induced matching of size $t+1$. We again distinguish cases according to the minimum degree of $G$.

\medskip

\textsc{Case I:} $G$ has a vertex $x$ of degree 1.

\medskip

Let $y$ be the neighbor of $x$. If $xy$ is an isolated edge, then $G\setminus \{x,y\}$ does not contain induced matchings of size $t$. Applying Observation~\ref{obs} and Fact~\ref{fact2} yields $\mis(G)\le 2\mis_{\bigtriangleup,t-1}(n-2) \leq  2 h_{t-1}(n-2)=h_t(n)$. If $y$ has further neighbors, then $G\setminus {N[y]}$ does not contain induced matchings of size $t$. Applying Observation~\ref{obs} and Fact~\ref{fact2} yields $\mis(G)\le \mis_{\bigtriangleup,t}(n-2) + \mis_{\bigtriangleup,t-1}(n-3) \leq  h_t(n-2)+h_{t-1}(n-3)\le h_t(n)$.

\medskip

\textsc{Case II:} $G$ has minimum degree $d \geq 4$.

\medskip

Applying Observation~\ref{obs} and Fact~\ref{fact2} yields $\mis(G)\le (d+1)h_{t}(n-d-1)\le h_t(n)$.

\medskip

\textsc{Case III:} $G$ has a component $C$ of minimum degree $3$ with a vertex of degree at least $4$.

\medskip

Let $x$ be a vertex of $C$ of degree $3$ and $y\in N(x)$ of degree at least $4$. Applying Observation~\ref{obs} and Fact~\ref{fact2} yields $\mis(G)\le 3h_t(n-4)+h_t(n-5)\le h_t(n)$.

\medskip

\textsc{Case IV:} $G$ is $3$-regular.

\medskip

Suppose first that there exist vertices $x,y$ with $N(x)=N(y)$. As $G$ is triangle-free, $x$ and $y$ cannot be adjacent. Moreover, for any maximal independent set $X$, we have $x\in X$ if and only if $y\in X$. Therefore, when applying Observation \ref{obs}, the number of maximal independent sets containing $x$ can be bounded by $\mis(G\setminus ({N[x]} \cup \{y\})$ instead of $\mis(G\setminus ({N[x]})$. We obtain $\mis(G)\le 3h_t(n-4)+h_t(n-5)\le h_t(n)$.

Suppose next that $G$ does not contain two vertices $x,y$ with $N(x)=N(y)$. Let $v$ be an arbitrary vertex of $G$ and $N(v)=\{a,b,c\}$. We apply Observation \ref{obs} in a slightly modified form: we keep $\mis(G\setminus {N[v]})$ and $\mis(G\setminus {N[a]})$, but to bound the number of maximal independent sets $X$ that contain $b$ or $c$, we use $\mis(G\setminus ({N[b]}\cup \{c\}))+\mis(G\setminus ({N[c]}\cup \{b\}))+\mis(G\setminus ({N[b]}\cup {N[c]}))$. The three terms bound the number of $X$s that contain exactly $b$, exactly $c$ or both $b$ and $c$, respectively. Observe that by the triangle-free property, $b$ and $c$ are not adjacent. Also, as $N(b)\neq N(c)$, we have $|{N[b]}\cup {N[c]}|\ge 6$. Therefore, we obtain $\mis(G)\le 2h_t(n-4)+2h_t(n-5)+h_t(n-6) \leq  h_t(n)$ by the definition of $h_t(n)$.

\medskip

\textsc{Case V}: $G$ has a component $C$ of minimum degree $2$ with a vertex of degree at least $3$.

\medskip

Let $x$ be a vertex of $C$ of degree $2$, and $y\in N(x)$ of degree at least $3$. Then applying Observation~\ref{obs} and Fact~\ref{fact2}, we obtain $\mis(G)\le 2h_t(n-3)+h_t(n-4)\le h_t(n)$.

\medskip

\textsc{Case VI}: $G$ is 2-regular, i.e., a cycle factor.

\medskip

As in Case IV of Theorem \ref{main}, we have $\mis(G)\le 5^{n/5}\le h_t(n)$.
\end{proof}


\begin{thebibliography}{99}


\bibitem{F}
Z. F\"uredi, The number of maximal independent sets in connected graphs. \textit{Journal of Graph Theory}, 11(4) (1987) 463-–470.


\bibitem{Gr} J.M. Griggs, C.M. Grinstead, D.R. Guichard, The number of maximal independent sets in a connected graph. \textit{Discrete Math.}, 68(2–3) (1988)
211--220.

\bibitem{HT}
M. Hujter, Zs. Tuza, The number of maximal independent sets in triangle-free graphs. \textit{SIAM Journal on Discrete Mathematics}, 6(2) (1993) 284--288.

\bibitem{KP}
J. Kahn, J. Park, Stability for Maximal Independent Sets. \textit{the Electronic Journal of Combinatorics}, 27(1) (2020) P1.59
\bibitem{MiMu}
R.E. Miller, D.E. Muller, A problem of maximum consistent subsets, IBM Research Report RC-240,
Thomas J. Watson Research Center, New York, USA, 1960.
\bibitem{MM}
J.W. Moon, L. Moser, On cliques in graphs. \textit{Israel Journal of Mathematics}, 3(1) (1965) 23--28.

\bibitem{Sa}
B. Sagan, A note on independent sets in trees. \textit{SIAM J. Discrete Math.}, 1(1) (1988) 105--108.


\bibitem{Wi} H.S. Wilf, The number of maximal independent sets in a tree. \textit{SIAM J. Algebraic Discrete Methods}, 7 (1986) 125--130.


\bibitem{W}
D.R Wood, On the number of maximal independent sets in a graph. \textit{Discrete Mathematics \& Theoretical Computer Science}, 13(3) (2011) 17--20.
\end{thebibliography}
\end{document}